\documentclass[11pt]{amsart}
\usepackage{amsmath,amssymb,mathrsfs}
\newtheorem{theorem}{Theorem}[section]
\newtheorem{corollary}[theorem]{Corollary}
\newtheorem{lemma}[theorem]{Lemma}
\newtheorem{proposition}[theorem]{Proposition}

\begin{document}

\title[Lagrangian minimal surface equation]{On the Lagrangian minimal surface equation and related problems}
\author{Simon Brendle}
\address{Department of Mathematics \\ Stanford University \\ 450 Serra Mall, Bldg 380 \\ Stanford, CA 94305} 
\thanks{The author was supported in part by the National Science Foundation under grant DMS-0905628.}
\begin{abstract}
We give a survey of various existence results for minimal Lagrangian graphs. We also discuss the mean curvature flow for Lagrangian graphs.
\end{abstract}
\maketitle

\section{Background on minimal Lagrangian geometry}

Minimal submanifolds are among the central objects in differential geometry. There is an important subclass of minimal submanifolds which was introduced by Harvey and Lawson \cite{Harvey-Lawson} in 1982. Given a Riemannian manifold $(M,g)$, a calibrating form $\Omega$ is a closed $m$-form on $M$ with the property that $\Omega(e_1,\hdots,e_m) \leq 1$ for each point $p \in M$ and every orthonormal $k$-frame $\{e_1,\hdots,e_m\} \subset T_p M$. An oriented $m$-dimensional submanifold $\Sigma \subset M$ is said to be calibrated by $\Omega$ if $\Omega(e_1,\hdots,e_m) = 1$ for every point $p \in \Sigma$ and every positively ortiented orthonormal basis $\{e_1,\hdots,e_m\}$ of $T_p \Sigma$. Using Stokes theorem, Harvey and Lawson showed that every calibrated submanifold is necessarily minimal:

\begin{theorem}[R.~Harvey, H.B.~Lawson \cite{Harvey-Lawson}]
Let $(M,g)$ be a Riemannian manifold. Moreover, let $\Omega$ be a calibrating $k$-form and let $\Sigma$ be a $k$-dimensional submanifold calibrated by $\Sigma$. Then $\Sigma$ minimizes volume in its homology class.
\end{theorem}

In the following, we consider the special case when $(M,g)$ is the Euclidean space $\mathbb{R}^{2n}$. We denote by $(x_1,\hdots,x_n,y_1,\hdots,y_n)$ the standard coordinates on $\mathbb{R}^{2n}$. Moreover, we denote by $\omega = \sum_{k=1}^n dx_k \wedge dy_k$ the standard symplectic form. Let $J$ be the associated complex structure, so that $J \frac{\partial}{\partial x_k} = \frac{\partial}{\partial y_k}$ and $J \frac{\partial}{\partial y_k} = -\frac{\partial}{\partial x_k}$. Finally, we define 
\[\sigma = (dx_1 + i \, dy_1) \wedge \hdots \wedge (dx_n + i \, dy_n).\] 
Note that $\sigma$ is a complex-valued $n$-form on $\mathbb{R}^{2n}$. Moreover, we have 
\[\sigma(Jv_1,v_2,\hdots,v_n) = i \, \sigma(v_1,v_2,\hdots,v_n)\] 
for all vectors $v_1,\hdots,v_n \in \mathbb{R}^n$.

Let now $\Sigma$ be a submanifold of $\mathbb{R}^{2n}$ of dimension $n$. Recall that $\Sigma$ is said to be Lagrangian if $\omega|_\Sigma = 0$. If $\Sigma$ is a Lagrangian submanifold, then it can be shown that $|\sigma(e_1,\hdots,e_n)| = 1$, where $\{e_1,\hdots,e_n\}$ is an orthonormal basis of $T_p \Sigma$. We may therefore write 
\begin{equation} 
\label{lagrangian.angle}
\sigma(e_1,\hdots,e_n) = e^{i\gamma} 
\end{equation}
for some function $\gamma: \Sigma \to \mathbb{R} / 2\pi\mathbb{Z}$. The function $\gamma$ is referred to as the Lagrangian angle of $\Sigma$.

The mean curvature vector of a Lagrangian submanifold $\Sigma$ is given by $J \, \nabla^\Sigma \gamma$, where $\nabla^\Sigma \gamma \in T_p\Sigma$ denotes the gradient of the Lagrangian angle. In particular, this implies:

\begin{theorem}[R.~Harvey, H.B.~Lawson \cite{Harvey-Lawson}]
If $\Sigma$ is a Lagrangian submanifold with $H = 0$, then the Lagrangian angle must be constant. Conversely, if $\Sigma$ is a Lagrangian and the Lagrangian angle is constant (so that $\gamma = c$), then $\Sigma$ is calibrated by the $n$-form $\Omega = \text{\rm Re}(e^{-ic} \, \sigma)$.
\end{theorem}

In particular, minimal Lagrangian submanifolds are special cases of calibrated submanifolds.

The first non-trivial examples of minimal Lagrangian submanifolds in $\mathbb{R}^{2n}$ were constructed by Harvey and Lawson \cite{Harvey-Lawson}. These examples are nearly flat and are constructed by means of the implicit function theorem.

\section{Minimal Lagrangian graphs in $\mathbb{R}^{2n}$}

We now assume that $\Sigma$ is an $n$-dimensional submanifold of $\mathbb{R}^{2n}$ which can be written as a graph over a Lagrangian plane in $\mathbb{R}^{2n}$. In other words, we write 
\[\Sigma = \{(x_1,\hdots,x_n,y_1,\hdots,y_n) \in \mathbb{R}^{2n}: (y_1,\hdots,y_n) = f(x_1,\hdots,x_n)\}.\] 
Here, the map $f$ is defined on some domain in $\mathbb{R}^n$ and takes values in $\mathbb{R}^n$. 

The condition that $\Sigma$ is Lagrangian is equivalent to the condition that $\partial_k f_l = \partial_l f_k$. Thus, $\Sigma$ is Lagrangian if and only if the map $f$ can locally be written as the gradient of some real-valued function $u$. In this case, the Lagrangian angle of $\Sigma$ is given by 
\[\gamma = \sum_{k=1}^n \arctan(\lambda_k),\] 
where $\lambda_1,\hdots,\lambda_k$ denote the eigenvalues of $Df(x) = D^2 u(x)$. Therefore, $\Sigma$ is a minimal Lagrangian submanifold if and only if $u$ satisfies the Hessian equation 
\begin{equation} 
\label{pde}
\sum_{k=1}^n \arctan(\lambda_k) = c. 
\end{equation}
A natural question is to classify all entire solutions of (\ref{pde}). In this direction Tsui and Wang proved the following result: 

\begin{theorem}[M.P.~Tsui, M.T.~Wang \cite{Tsui-Wang1}]
Let $f: \mathbb{R}^n \to \mathbb{R}^n$ be a smooth map with the property that the graph $\Sigma = \{(x,f(x)): x \in \mathbb{R}^n\}$ is a minimal Lagrangian submanifold. Moreover, we assume that, for each point $x \in \mathbb{R}^n$, the eigenvalues of $Df(x)$ satisfy $\lambda_i \lambda_j \geq -1$ and $|\lambda_i| \leq K$. Then $f$ is an affine function. 
\end{theorem}

A closely related Bernstein-type result was established independently in \cite{Yuan}:

\begin{theorem}[Y.~Yuan \cite{Yuan}]
Let $u: \mathbb{R}^n \to \mathbb{R}$ be a smooth convex solution of (\ref{pde}). Then $u$ is a quadratic polynomial.
\end{theorem}

In order to study the equation (\ref{pde}) on a bounded domain in $\mathbb{R}^n$, one needs to impose a boundary condition. One possibility is to impose a Dirichlet boundary condition for the potential function $u$. This boundary value problem was studied in the fundamental work of Caffarelli, Nirenberg, and Spruck \cite{Caffarelli-Nirenberg-Spruck}. In particular, they obtained the following existence theorem: 

\begin{theorem}[L.~Caffarelli, L.~Nirenberg, J.~Spruck \cite{Caffarelli-Nirenberg-Spruck}]
Let $\Omega$ be a uniformly convex domain in $\mathbb{R}^n$, and let $\varphi: \partial \Omega \to \mathbb{R}$ be a smooth function. Then there exists a smooth function $u: \Omega \to \mathbb{R}$ satisfying 
\[\sum_{k=1}^n \arctan(\lambda_k) = \Big [ \frac{n-1}{2} \Big ] \, \pi\] 
and $u|_{\partial \Omega} = \varphi$.
\end{theorem}

We now describe another natural boundary condition for (\ref{pde}). Instead of prescribing the boundary values of $u$, we prescribe the image of $\Omega$ under the map $f = \nabla u$. This choice of boundary condition has been studied before in connection with the Monge-Amp\`ere equation and other Hessian equations (see \cite{Caffarelli}, \cite{Urbas1}, \cite{Urbas2}).

\begin{theorem}[S.~Brendle, M.~Warren \cite{Brendle-Warren}] 
\label{brendle.warren.1}
Let $\Omega$ and $\tilde{\Omega}$ be uniformly convex domains in $\mathbb{R}^n$. Then we can find a smooth function $u: \Omega \to \mathbb{R}$ and a real number $c$ with the following properties: 
\begin{itemize}
\item[(i)] The function $u$ is uniformly convex.
\item[(ii)] The function $u$ solves the equation (\ref{pde}).
\item[(iii)] The map $\nabla u: \Omega \to \mathbb{R}$ is a diffeomorphism from $\Omega$ to $\tilde{\Omega}$.
\end{itemize}
Moreover, the pair $(u,c)$ is unique.
\end{theorem}

Thus, we can draw the following conclusion:

\begin{corollary}[S.~Brendle, M.~Warren \cite{Brendle-Warren}] 
\label{brendle.warren.2}
Let $\Omega$ and $\tilde{\Omega}$ be uniformly convex domains in $\mathbb{R}^n$ with smooth boundary. Then there exists a diffeomorphism $f: \Omega \to \tilde{\Omega}$ such that the graph $\Sigma = \{(x,f(x)): x \in \Omega\}$ is a minimal Lagrangian submanifold of $\mathbb{R}^{2n}$.
\end{corollary}

In particular, the submanifold $\Sigma$ satisfies $\partial \Sigma \subset \partial \Omega \times \partial \tilde{\Omega}$. Thus, the surface $\Sigma$ satisfies a free boundary value problem.

We note that the potential function $u$ is not a geometric quantity; on the other hand, the gradient $\nabla u = f$ does have geometric significance. From a geometric point of view, the second boundary value problem is more natural than the Dirichlet boundary condition.

We now describe the proof of Theorem \ref{brendle.warren.1}. The uniqueness statement follows from a standard argument based on the maximum principle. In order to prove the existence statement, we use the continuity method. The idea is to deform $\Omega$ and $\tilde{\Omega}$ to the unit ball in $\mathbb{R}^n$. As usual, the central issue is to bound the Hessian of the potential function $u$. In geometric terms, this corresponds to a bound on the slope of $\Sigma$. 

\begin{proposition}[\cite{Brendle-Warren}]
\label{a.priori.estimates}
Let us fix two uniformly convex domains $\Omega$ and $\tilde{\Omega}$. Moreover, let $u$ be a convex solution of (\ref{pde}) with the property that $\nabla u$ is a diffeomorphism from $\Omega$ to $\tilde{\Omega}$. Then $|D^2 u(x)| \leq C$ for all points $x \in \Omega$ and all vectors $v \in \mathbb{R}^n$. Here, $C$ is a positive constant, which depends only on $\Omega$ and $\tilde{\Omega}$.
\end{proposition}

The proof of Proposition \ref{a.priori.estimates} is inspired by earlier work of Urbas \cite{Urbas1}, \cite{Urbas2}. By assumption, we can find uniformly convex boundary defining functions $h: \Omega \to (-\infty,0]$ and $\tilde{h}: \tilde{\Omega} \to (-\infty,0]$, so that $h|_{\partial \Omega} = 0$ and $\tilde{h}|_{\partial \tilde{\Omega}} = 0$. Moreover, let us fix a constant $\theta > 0$ such that $D^2 h(x) \geq \theta \, I$ for all points $x \in \Omega$ and $D^2 \tilde{h}(y) \geq \theta \, I$ for all points $y \in \tilde{\Omega}$. 

In the following, we sketch the main steps involved in the proof of Proposition \ref{a.priori.estimates}. 

\textit{Step 1:} Let $u$ be a convex solution of (\ref{pde}) with the property that $\nabla u$ is a diffeomorphism from $\Omega$ to $\tilde{\Omega}$. Differentiating the equation (\ref{pde}), we obtain 
\begin{equation} 
\label{pde.2}
\sum_{i,j=1}^n a_{ij}(x) \, \partial_i \partial_j \partial_k u(x) = 0 
\end{equation}
for all $x \in \Omega$ and all $k \in \{1,\hdots,n\}$. Here, the coefficients $a_{ij}(x)$ are defined as the components of the matrix $A(x) = (I + (D^2 u(x))^2)^{-1}$.

We now define a function $H: \Omega \to \mathbb{R}$ by $H(x) = \tilde{h}(\nabla u(x))$. Using the identity (\ref{pde.2}), one can show that 
\[\bigg | \sum_{i,j=1}^n a_{ij}(x) \, \partial_i \partial_j H(x) \bigg | \leq C\] 
for some uniform constant $C$. Using the maximum principle, we conclude that $H(x) \geq C \, h(x)$ for all points $x \in \Omega$. Here, $C$ is a uniform constant which depends only on $\Omega$ and $\tilde{\Omega}$. This implies $\langle \nabla h(x),\nabla H(x) \rangle \leq C \, |\nabla h(x)|^2$ at each point $x \in \partial \Omega$. As a result, we can bound certain components of the Hessian of $u$ along $\partial \Omega$.

\textit{Step 2:} In the next step, we prove a uniform obliqueness estimate. To that end, we consider the function $\chi(x) = \langle \nabla h(x),\nabla \tilde{h}(\nabla u(x)) \rangle$. It is not difficult to show that $\chi(x) > 0$ for all $x \in \partial \Omega$. The goal is to obtain a uniform lower bound for $\inf_{x \in \partial \Omega} \chi(x)$. Using the relation (\ref{pde.2}), one can show that 
\[\bigg | \sum_{i,j=1}^n a_{ij}(x) \, \partial_i \partial_j \chi(x) \bigg | \leq C\] 
for some uniform constant $C$. We can therefore find a uniform constant $K$ such that 
\[\sum_{i,j=1}^n a_{ij}(x) \, \partial_i \partial_j (\chi(x) - K \, h(x)) \leq 0.\] 
We now consider a point $x_0 \in \partial \Omega$, where the function $\chi(x) - K \, h(x)$ attains its global minimum. Then $\nabla \chi(x_0) = (K-\mu) \, \nabla h(x_0)$ for some real number $\mu \geq 0$. Hence, we obtain 
\begin{align*} 
(K-\mu) \, \chi(x_0) 
&= \langle \nabla \chi(x_0),\nabla \tilde{h}(\nabla u(x_0)) \rangle \\ 
&= \sum_{i,j=1}^n \partial_i \partial_j h(x_0) \, (\partial_i \tilde{h})(\nabla u(x_0)) \, (\partial_j \tilde{h})(\nabla u(x_0)) \\ 
&+ \sum_{i,j=1}^n (\partial_i \partial_j \tilde{h})(\nabla u(x_0)) \, \partial_i h(x_0) \, \partial_j H(x_0) \\ 
&\geq \theta \, |\nabla \tilde{h}(\nabla u(x_0))|^2 + \sum_{i,j=1}^n (\partial_i \partial_j \tilde{h})(\nabla u(x_0)) \, \partial_i h(x_0) \, \partial_j H(x_0).
\end{align*} 
Since $\nabla H(x_0)$ is a positive multiple of $\nabla h(x_0)$, it follows that 
\[K \, \chi(x_0) \geq \theta \, |\nabla \tilde{h}(\nabla u(x_0))|^2.\] 
Since $\inf_{x \in \partial \Omega} \chi(x) = \chi(x_0)$, we obtain a uniform lower bound for $\inf_{x \in \partial \Omega} \chi(x)$.

\textit{Step 3:} Having established the uniform obliqueness estimate, we next bound the tangential components of the Hessian $D^2 u(x)$ for each point $x \in \partial \Omega$. To explain this, let 
\[M = \sup \bigg \{ \sum_{k,l=1}^n \partial_k \partial_l u(x) \, z_k \, z_l: x \in \partial \Omega, \, z \in T_x (\partial \Omega), \, |z| = 1 \bigg \}.\] 
Our goal is to establish an upper bound for $M$. To that end, we fix a point $x_0 \in \partial M$ and a vector $w \in T_{x_0}(\partial \Omega)$ such that $|w| = 1$ and 
\[\sum_{k,l=1}^n \partial_k \partial_l u(x_0) \, w_k \, w_l = M.\] 
We then consider the function 
\[\psi(x) = \sum_{k,l=1}^n \partial_k \partial_l u(x) \, w_k \, w_l.\] 
Differentiating the identity (\ref{pde}) twice, we obtain 
\[\sum_{i,j=1}^n a_{ij}(x) \, \partial_i \partial_j \psi(x) \geq 0\] 
for all $x \in \Omega$. Using the definition of $M$, it can be shown that 
\begin{align} 
\label{boundary.estimate}
\psi(x) &\leq M \, \bigg | w - \frac{\langle \nabla h(x),w \rangle}{\langle \nabla h(x),\nabla \tilde{h}(\nabla u(x)) \rangle} \, \nabla \tilde{h}(\nabla u(x)) \bigg |^2 \notag\\ 
&+ L \, \langle \nabla h(x),w \rangle^2 
\end{align} 
for all points $x \in \partial \Omega$. Here, $L$ is fixed constant that depends only on $\Omega$ and $\tilde{\Omega}$. 

Let $\varepsilon$ be a positive real number such that $\inf_{x \in \partial \Omega} \chi(x) > \varepsilon$, and let $\eta: \mathbb{R} \to (0,\infty)$ be a smooth function satisfying $\eta(s) = s$ for all $s \geq \varepsilon$. Using (\ref{boundary.estimate}) and the maximum principle, we obtain an estimate of the form 
\begin{align} 
\label{interior.estimate}
\psi(x) &\leq M \, \bigg | w - \frac{\langle \nabla h(x),w \rangle}{\eta(\langle \nabla h(x),\nabla \tilde{h}(\nabla u(x)) \rangle)} \, \nabla \tilde{h}(\nabla u(x)) \bigg |^2 \notag \\ 
&+ L \, \langle \nabla h(x),w \rangle^2 - C \, h(x)
\end{align} 
for all $x \in \Omega$. Moreover, equality holds in (\ref{interior.estimate}) when $x = x_0$. Consequently, we obtain a lower bound for the normal derivative of $\psi$ at the point $x_0$. More precisely, 
\[\langle \nabla \psi(x_0),\nabla \tilde{h}(\nabla u(x_0)) \rangle + C \, M + C \geq 0,\] 
where $C$ is a uniform constant that depends only on $\Omega$ and $\tilde{\Omega}$. On the other hand, we have 
\begin{align*} 
&\langle \nabla \psi(x_0),\nabla \tilde{h}(\nabla u(x_0)) \rangle + \theta \, M^2 \\ 
&\leq \sum_{i,k,l=1}^n (\partial_i \tilde{h})(\nabla u(x_0)) \, \partial_i \partial_k \partial_l u(x_0) \, w_k \, w_l \\ 
&+ \sum_{i,j,k,l=1}^n (\partial_i \partial_j \tilde{h})(\nabla u(x_0)) \, \partial_i \partial_k u(x_0) \, \partial_j \partial_l u(x_0) \, w_k \, w_l \\ 
&= \sum_{k,l=1}^n \partial_k \partial_l H(x_0) \, w_k \, w_l \\ 
&= -\langle \nabla H(x_0),I\!I(w,w) \rangle, 
\end{align*} 
where $I\!I$ denotes the second fundamental form of $\partial \Omega$. Consequently, 
\[\langle \nabla \psi(x_0),\nabla \tilde{h}(\nabla u(x_0)) \rangle + \theta \, M^2 \leq C.\] 
Putting these facts together, we obtain an a-priori estimate for $M$.

\textit{Step 4:} Once we have uniform bounds for the Hessian of $u$ along the boundary, we can use the maximum principle to bound the Hessian of $u$ in the interior of $\Omega$. This step is by now standard, and follows ideas in \cite{Caffarelli-Nirenberg-Spruck}.

\section{Area-preserving minimal maps between surfaces}

We now describe a different boundary problem value for minimal Lagrangian graphs. To that end, let $M$ be a two-dimensional surface equipped with a Riemannian metric $g$ and a complex structure $J$. We consider the product $M = N \times N$ equipped with the product metric. We define a complex structure on $M$ by 
\[J_{(p,q)}(w,\tilde{w}) = (J_p w,-J_q \tilde{w})\] 
for all vectors $w \in T_p N$ and $\tilde{w} \in T_q N$. 

Our goal is to construct minimal Lagrangian submanifolds in $M$. We will assume throughout this section that $N$ is a surface with constant Gaussian curvature, so that $M$ is a K\"ahler-Einstein manifold. (Otherwise, the minimal Lagrangian equation leads to an overdetermined system of PDEs). 

In the special case when $N = \mathbb{R}^2$, the existence of minimal Lagrangian graphs can be reduced to the solvability of the second boundary value problem for the Monge-Amp\`ere equation. To describe this, we consider two domains $\Omega,\tilde{\Omega} \subset \mathbb{R}^2$. Moreover, we consider a diffeomorphism $f: \Omega \to \tilde{\Omega}$, and let 
\[\Sigma = \{(p,f(p)): p \in \Omega\}.\] 
The graph $\Sigma$ is Lagrangian if and only if the map $f$ is area-preserving and orientation-preserveing, so that $\det Df = 1$. Moreover, $\Sigma$ has vanishing mean curvature if and only if the Lagrangian angle is constant; this means that 
\[\cos \gamma \, (\partial_1 f_2 - \partial_2 f_1) = \sin \gamma \, (\partial_1 f_1 + \partial_2 f_2)\] 
for some constant $\gamma \in \mathbb{R}$. Hence, we may locally write 
\begin{align*} 
f_1 &= \cos \gamma \, \partial_1 u - \sin \gamma \, \partial_2 u \\ 
f_2 &= \sin \gamma \, \partial_1 u + \cos \gamma \, \partial_2 u 
\end{align*} 
for some potential function $u$.

In other words, the map $f$ can locally be expressed as the composition of a gradient mapping with a rotation in $\mathbb{R}^2$. Since $f$ is area-preserving, the potential function solves the Monge-Amp\`ere equation $\det D^2 u = 1$.

It was shown by Delano\"e \cite{Delanoe} that the second boundary value problem for the Monge-Amp\`ere equation is solvable, provided that $\Omega$ and $\tilde{\Omega}$ are uniformly convex and have the same area. This implies the following result:

\begin{theorem}[P.~Delano\"e \cite{Delanoe}] 
Let $\Omega$ and $\tilde{\Omega}$ be uniformly convex domains in $\mathbb{R}^2$ with smooth boundary. Assume that $\Omega$ and $\tilde{\Omega}$ have the same area. Then there exists a minimal Lagrangian diffeomorphism from $\Omega$ to $\tilde{\Omega}$.
\end{theorem}

The assumption that $\Omega$ and $\tilde{\Omega}$ are uniformly convex cannot be removed. In fact, Urbas \cite{Urbas3} constructed two domains in $\mathbb{R}^2$ such that the second boundary value for the Monge-Amp\`ere equation does not admit a smooth solution. In this example, the domain $\Omega$ is the unit disk; moreover, the geodesic curvature of $\partial \tilde{\Omega}$ is greater than $-\varepsilon$.

We next consider the case when $N$ is a complete, simply connected surface with negative Gaussian curvature. In this case, we have the following result:

\begin{theorem}[S.~Brendle \cite{Brendle}]
\label{existence.uniqueness}
Let $N$ be a complete, simply connected surface with constant negative Gaussian curvature, and let $\Omega$ and $\tilde{\Omega}$ be uniformly convex domains in $N$ with smooth boundary. Assume that $\Omega$ and $\tilde{\Omega}$ have the same area. Given any point $\overline{p} \in \partial \Omega$ and any point $\overline{q} \in \partial \tilde{\Omega}$, there exists a unique minimal Lagrangian diffeomorphism from $\Omega$ to $\tilde{\Omega}$ that maps $\overline{p}$ to $\overline{q}$.
\end{theorem}

We note that the product $M$ does not admit a parallel complex volume form. Therefore, we do not have a notion of Lagrangian angle in this setting. As a result, it is no longer possible to reduce the minimal Lagrangian equation to a PDE for a scalar function.

The proof of Theorem \ref{existence.uniqueness} uses the continuity method. In order to make the continuity argument work, it is necessary to establish a-priori estimates for area-preserving minimal maps between domains in $N$. Among other things, the proof uses the following lemma, which was first obtained by Wang \cite{Wang2} in his study of the Lagrangian mean curvature flow:

\begin{lemma}[M.T.~Wang \cite{Wang2}]
\label{pde.for.beta}
Suppose that $f: \Omega \to \tilde{\Omega}$ is an area-preserving minimal map. Moreover, let $\Sigma = \{(p,f(p)): p \in \Omega\}$ denote the graph of $f$, and let $\beta: \Sigma \to \mathbb{R}$ by defined by 
\[\beta(p,f(p)) = \frac{2}{\sqrt{\det(I + Df_p^* \, Df_p)}}.\] 
Then 
\[\Delta_\Sigma \beta = -2 \, |I\!I|^2 \, \beta - \kappa \, \beta \, (1 - \beta^2).\] 
Here, $\kappa < 0$ denotes the Gaussian curvature of the two-dimensional surface $N$. 
\end{lemma}

We now describe the proof of Lemma \ref{pde.for.beta}. Given any point $(p,q) \in M$, we define a two-form $\rho: T_{(p,q)} M \times T_{(p,q)} M \to \mathbb{R}$ by 
\[\rho \big ( (w_1,\tilde{w}_1),(w_2,\tilde{w}_2) \big ) = \langle Jw_1,w_2 \rangle + \langle J\tilde{w}_1,\tilde{w}_2 \rangle\] 
for all vectors $w_1,w_2 \in T_p N$ and $\tilde{w}_1,\tilde{w}_2 \in T_q N$. Clearly, $\rho$ is parallel. At each point on $\Sigma$, we have $\beta = \rho(e_1,e_2)$, where $\{e_1,e_2\}$ is a local orthonormal frame for $T\Sigma$. Differentiating this identity, we obtain 
\[V(\beta) = \rho(I\!I(e_1,V),e_2) + \rho(e_1,I\!I(e_2,V))\]
for every vector $V \in T\Sigma$. This implies 
\begin{align} 
\label{laplacian.beta}
\Delta_\Sigma \beta 
&= \sum_{k=1}^2 \rho(\nabla_{e_k}^M I\!I(e_1,e_k),e_2) + \sum_{k=1}^2 \rho(e_1,\nabla_{e_k}^M I\!I(e_2,e_k)) \notag \\ 
&+ 2 \sum_{k=1}^2 \rho(I\!I(e_1,e_k),I\!I(e_2,e_k)). 
\end{align} 
Using the Codazzi equations (see e.g. \cite{ONeill}, Chapter 4, Proposition 33) we obtain 
\begin{align} 
\label{aux.1}
&\sum_{k=1}^2 \rho(\nabla_{e_k}^M I\!I(e_1,e_k),e_2) + \sum_{k=1}^2 \rho(e_1,\nabla_{e_k}^M I\!I(e_2,e_k)) \notag \\ 
&= \sum_{k=1}^2 \rho(\nabla_{e_k}^\perp I\!I(e_1,e_k),e_2) + \sum_{k=1}^2 \rho(e_1,\nabla_{e_k}^\perp I\!I(e_2,e_k)) \notag \\ 
&+ \sum_{k=1}^2 \langle \nabla_{e_k}^M I\!I(e_1,e_k),e_1 \rangle \, \rho(e_1,e_2) + \sum_{k=1}^2 \langle \nabla_{e_k}^M I\!I(e_2,e_k),e_2 \rangle \, \rho(e_1,e_2) \notag \\ 
&= \sum_{k=1}^2 \rho(\nabla_{e_1}^\perp I\!I(e_k,e_k),e_2) + \sum_{k=1}^2 \rho(e_1,\nabla_{e_2}^\perp I\!I(e_k,e_k)) \\ 
&- \sum_{k=1}^2 |I\!I(e_1,e_k)|^2 \, \rho(e_1,e_2) - \sum_{k=1}^2 |I\!I(e_2,e_k)|^2 \, \rho(e_1,e_2) \notag \\ 
&- R_M(e_2,e_1,e_2,Je_1) \, \rho(Je_1,e_2) - R_M(e_2,e_1,e_2,Je_2) \, \rho(Je_2,e_2) \notag \\ 
&- R_M(e_1,e_2,e_1,Je_1) \, \rho(e_1,Je_1) - R_M(e_1,e_2,e_1,Je_2) \, \rho(e_1,Je_2). \notag
\end{align} 
Here, $\nabla^\perp$ denotes the induced connection on the normal bundle of $\Sigma$. Since $N$ has constant Gaussian curvature $\kappa$, we have 
\begin{align*} 
&R_M(e_2,e_1,e_2,Je_1) \, \rho(Je_1,e_2) + R_M(e_2,e_1,e_2,Je_2) \, \rho(Je_2,e_2) \\ 
&+ R_M(e_1,e_2,e_1,Je_1) \, \rho(e_1,Je_1) + R_M(e_1,e_2,e_1,Je_2) \, \rho(e_1,Je_2) \\ 
&= \kappa \, \beta \, (1-\beta^2). 
\end{align*} 
Substituting this into (\ref{aux.1}) gives 
\begin{align} 
\label{aux.2}
&\sum_{k=1}^2 \rho(\nabla_{e_k}^M I\!I(e_1,e_k),e_2) + \sum_{k=1}^2 \rho(e_1,\nabla_{e_k}^M I\!I(e_2,e_k)) \notag \\ 
&= -|I\!I|^2 \, \beta - \kappa \, \beta \, (1-\beta^2). 
\end{align}
Moreover, we have 
\begin{align} 
\label{aux.3}
&\sum_{k=1}^2 \rho(I\!I(e_1,e_k),I\!I(e_2,e_k)) \notag \\ 
&= \sum_{k=1}^2 \langle I\!I(e_1,e_k),Je_1 \rangle \, \langle I\!I(e_2,e_k),Je_2 \rangle \, \rho(Je_1,Je_2) \notag \\ 
&+ \sum_{k=1}^2 \langle I\!I(e_1,e_k),Je_2 \rangle \, \langle I\!I(e_2,e_k),Je_1 \rangle \, \rho(Je_2,Je_1) \notag \\ 
&= \sum_{k=1}^2 \langle I\!I(e_1,e_1),Je_k \rangle \, \langle I\!I(e_2,e_2),Je_k \rangle \, \beta \\ 
&- \sum_{k=1}^2 \langle I\!I(e_1,e_2),Je_k \rangle \, \langle I\!I(e_1,e_2),Je_k \rangle \, \beta \notag \\ 
&= -\frac{1}{2} \, |I\!I|^2 \, \beta. \notag
\end{align}
Combining (\ref{laplacian.beta}), (\ref{aux.2}), and (\ref{aux.3}), we conclude that 
\[\Delta_\Sigma \beta = -2 \, |I\!I|^2 \, \beta - \kappa \, \beta \, (1 - \beta^2),\] 
as claimed. This completes the proof of Lemma \ref{pde.for.beta}. \\

We next describe the a key priori estimate for the differential $Df$.

\begin{proposition}[\cite{Brendle}]
\label{estimate.for.Df}
Let $\Omega$ and $\tilde{\Omega}$ be uniformly convex domains in $N$ with smooth boundary. Suppose that $f: \Omega \to \tilde{\Omega}$ is an area-preserving minimal map. Then $|Df_p| \leq C$ for all points $p \in \Omega$, where $C$ is a uniform constant that depends only on $\Omega$ and $\tilde{\Omega}$.
\end{proposition}

Let us sketch the main ideas involved in the proof of Proposition \ref{estimate.for.Df}. Let $h: \Omega \to (-\infty,0]$ and $\tilde{h}: \tilde{\Omega} \to (-\infty,0]$ be uniformly convex boundary defining functions for $\Omega$ and $\tilde{\Omega}$. We may choose $h$ and $\tilde{h}$ such that $|\nabla h_p| = 1$ for all $p \in \partial \Omega$ and $|\nabla \tilde{h}_q| = 1$ for all $q \in \partial \tilde{\Omega}$.

Since $h$ and $\tilde{h}$ are uniformly convex, we have 
\begin{equation} 
\label{convexity.h}
\theta \, g \leq D^2 h \leq \frac{1}{\theta} \, g 
\end{equation}
and 
\begin{equation} 
\label{convexity.tilde.h}
\theta \, g \leq D^2 \tilde{h} \leq \frac{1}{\theta} \, g 
\end{equation}
for some positive constant $\theta$. 

\textit{Step 1:} Let $\Sigma = \{(p,f(p)): p \in \Omega\}$ denote the graph of $f$. By assumption, $\Sigma$ is a minimal submanifold of $M$. We next define two functions $H,\tilde{H}: \Sigma \to \mathbb{R}$ by $H(p,f(p)) = h(p)$ and $\tilde{H}(p,f(p)) = \tilde{h}(f(p))$. The relations (\ref{convexity.h}) and (\ref{convexity.tilde.h}) imply $\theta \leq \Delta_\Sigma H \leq \frac{1}{\theta}$ and $\theta \leq \Delta_\Sigma \tilde{H} \leq \frac{1}{\theta}$. Using the maximum principle, we obtain $\frac{1}{\theta^2} \, H \leq \tilde{H} \leq \theta^2 \, H$ at each point on $\Sigma$. In other words, we have 
\[\frac{1}{\theta^2} \, h(p) \leq \tilde{h}(f(p)) \leq \theta^2 \, h(p)\] 
for all points $p \in \Omega$. Consequently, 
\[\theta^2 \leq \langle Df_p(\nabla h_p),\nabla \tilde{h}_{f(p)} \rangle \leq \frac{1}{\theta^2}\]  for all points $p \in \partial \Omega$. 

\textit{Step 2:} In the next step, we define a linear isometry $Q_p: T_p N \to T_{f(p)} N$ by
\[Q_p = Df_p \, \big [ Df_p^* \, Df_p \big ]^{-\frac{1}{2}}.\] 
It is straightforward to verify that $J_{f(p)} \, Q_p = Q_p \, J_p$ for all $p \in \Omega$. We next define a bilinear form $\sigma: T_{(p,f(p))} M \times T_{(p,f(p))} M \to \mathbb{C}$ by 
\begin{align*} 
\sigma \big ( (w_1,\tilde{w}_1),(w_2,\tilde{w}_2) \big ) 
&= i \, \langle Q_p(w_1),\tilde{w}_2 \rangle + \langle Q_p(J_p w_1),\tilde{w}_2 \rangle \\ 
&- i \, \langle Q_p(w_2),\tilde{w}_1 \rangle - \langle Q_p(J_p w_2),\tilde{w}_1 \rangle 
\end{align*} 
for all vectors $w_1,w_2 \in T_p N$ and all vectors $\tilde{w}_1,\tilde{w}_2 \in T_{f(p)} N$. The bilinear form $\sigma$ satisfies $\sigma(W_2,W_1) = -\sigma(W_1,W_2)$ and $\sigma(JW_1,W_2) = i \, \sigma(W_1,W_2)$ for all vectors $W_1,W_2 \in T_{(p,f(p))} M$. Moreover, if $\{e_1,e_2\}$ is an orthonormal basis of $T_{(p,f(p))} \Sigma$, then $\sigma(e_1,e_2) = \pm 1$.

The crucial observation is that $\sigma$ is parallel with respect to the Levi-Civita connection on $M$. More precisely, suppose that $W_1$ and $W_2$ are vector fields on $M$. Then the expression $\sigma(W_1,W_2)$ defines a complex-valued function on $\Sigma$. The derivative of that function is given by 
\begin{equation} 
\label{parallel}
V(\sigma(W_1,W_2)) = \sigma(\nabla_V^M W_1,W_2) + \sigma(W_1,\nabla_V^M W_2). 
\end{equation}
The relation (\ref{parallel}) is a consequence of the fact that $\Sigma$ has zero mean curvature (see \cite{Brendle}, Proposition 3.3, for details). Differentiating the identity (\ref{parallel}), we obtain 
\begin{align} 
\label{laplacian} 
\Delta_\Sigma(\sigma(W_1,W_2)) 
&= \sum_{k=1}^2 \sigma(\nabla_{e_k,e_k}^{M,2} W_1,W_2) \notag \\ 
&+ \sum_{k=1}^2 \sigma(W_1,\nabla_{e_k,e_k}^{M,2} W_2) \\ 
&+ 2 \, \sum_{k=1}^2 \sigma(\nabla_{e_k}^M W_1,\nabla_{e_k}^M W_2). \notag
\end{align} 

\textit{Step 3:} We now define a function $\varphi: \Sigma \to \mathbb{R}$ by 
\[\varphi(p,f(p)) = \langle Q_p(\nabla h_p),\nabla \tilde{h}_{f(p)} \rangle.\] 
It is easy to see that $\varphi(p,f(p)) > 0$ for $p \in \partial \Omega$. Our goal is to establish a lower bound for $\inf_{p \in \partial \Omega} \varphi(p,f(p))$. This estimate can be viewed as a generalization of the uniform obliqueness estimate in \cite{Delanoe}.

To prove this estimate, we define vector fields $W_1$ and $W_2$ on $M$ by $(W_1)_{(p,q)} = (\nabla h_p,0)$ and $(W_2)_{(p,q)} = (0,\nabla \tilde{h}_q)$. Clearly, $\varphi = \text{\rm Re}(\sigma(W_1,W_2))$. Hence, the identity (\ref{laplacian}) implies $\Delta_\Sigma \varphi \leq L$, 
where $L$ is a positive constant that depends only on $\Omega$ and $\tilde{\Omega}$. Hence, we obtain $\Delta_\Sigma \big ( \varphi - \frac{L}{\theta} \, H \big ) \leq 0$. Consequently, the function $\varphi - \frac{L}{\theta} \, H$ attains its maximum at some point $(p_0,f(p_0)) \in \partial \Sigma$. At the point $(p_0,f(p_0))$, we have 
\[\nabla^\Sigma \varphi = \Big ( \frac{L}{\theta} - \mu \Big ) \, \nabla^\Sigma H\] 
for some real number $\mu \geq 0$. Consequently, for every vector $v \in T_{p_0} N$, we have 
\begin{align*} 
\Big ( \frac{L}{\theta} - \mu \Big ) \, \langle \nabla h_{p_0},v \rangle 
&= \Big ( \frac{L}{\theta} - \mu \Big ) \, \langle \nabla^\Sigma H,(v,Df_{p_0}(v)) \rangle \\ 
&= \big \langle \nabla^\Sigma \varphi,(v,Df_{p_0}(v)) \big \rangle \\ 
&= (D^2 h)_{p_0} \big ( v,Q_{p_0}^*(\nabla \tilde{h}_{f(p_0)}) \big ) \\ 
&+ (D^2 \tilde{h})_{f(p_0)} \big ( Q_{p_0}(\nabla h_{p_0}),Df_{p_0}(v) \big ). 
\end{align*} 
In particular, if we choose $v = Q_{p_0}^*(\nabla \tilde{h}_{f(p_0)})$, then we obtain 
\begin{align*} 
\Big ( \frac{L}{\theta} - \mu \Big ) \, \varphi(p_0,f(p_0)) 
&= (D^2 h)_{p_0} \big ( Q_{p_0}^*(\nabla \tilde{h}_{f(p_0)}),Q_{p_0}^*(\nabla \tilde{h}_{f(p_0)}) \big ) \\ 
&+ (D^2 \tilde{h})_{f(p_0)} \big ( Q_{p_0}(\nabla h_{p_0}),Q_{p_0}(Df_{p_0}^*(\nabla \tilde{h}_{f(p_0)})) \big ). 
\end{align*} 
By (\ref{convexity.h}), we have 
\begin{align*} 
&(D^2 h)_{p_0} \big ( Q_{p_0}^*(\nabla \tilde{h}_{f(p_0)}),Q_{p_0}^*(\nabla \tilde{h}_{f(p_0)}) \big ) \\ 
&\geq \theta \, |Q_{p_0}^*(\nabla \tilde{h}_{f(p_0)})|^2 = \theta \, |\nabla \tilde{h}_{f(p_0)}|^2 = \theta. 
\end{align*}
Moreover, the vector $Df_{p_0}^*(\nabla \tilde{h}_{f(p_0)})$ is a positive multiple of $\nabla h_{p_0}$. Since $\tilde{h}$ is convex, it follows that 
\[(D^2 \tilde{h})_{f(p_0)} \big ( Q_{p_0}(\nabla h_{p_0}),Q_{p_0}(Df_{p_0}^*(\nabla \tilde{h}_{f(p_0)})) \big ) \geq 0.\] 
Putting these facts together yields 
\[\Big ( \frac{L}{\theta} - \mu \Big ) \, \varphi(p_0,f(p_0)) \geq \theta,\] 
hence 
\begin{equation} 
\label{obliqueness}
\inf_{p \in \partial \Omega} \varphi(p,f(p)) = \varphi(p_0,f(p_0)) \geq \frac{\theta^2}{L}. 
\end{equation}

\textit{Step 4:} We next show that $|Df_p| \leq C$ for all points $p \in \partial \Omega$. To see this, let us define $v_1 = \nabla h_p$ and $v_2 = J \, \nabla h_p$. Similarly, we define $\tilde{v}_1 = \nabla \tilde{h}_{f(p)}$ and $\tilde{v}_2 = J \, \nabla \tilde{h}_{f(p)}$. Clearly, the vectors $\{v_1,v_2\}$ form an orthonormal basis of $T_p N$, and the vectors $\{\tilde{v}_1,\tilde{v}_2\}$ form an orthonormal basis of $T_{f(p)} N$. We now write 
\[Df_p(v_1) = a \, \tilde{v}_1 + b \, \tilde{v}_2\] 
and 
\[Df_p(v_2) = c \, \tilde{v}_2\] 
for suitable coefficients $a,b,c$. Note that $ac = 1$ since $f$ is area-preserving. Using the inequality $\theta^2 \leq \langle Df_p(\nabla h_p),\nabla \tilde{h}_{f(p)} \rangle \leq \frac{1}{\theta}^2$, we conclude that $\theta^2 \leq a \leq \frac{1}{\theta^2}$ and $\theta^2 \leq c \leq \frac{1}{\theta^2}$. In order to bound $b$, we observe that 
\begin{align*} 
a \, \langle Q_p(v_2),\tilde{v}_1 \rangle + b \, \langle Q_p(v_2),\tilde{v}_2 \rangle 
&= \langle Q_p(v_2),Df_p(v_1) \rangle \\ 
&= \langle Q_p(v_1),Df_p(v_2) \rangle \\ 
&= c \, \langle Q_p(v_1),\tilde{v}_2 \rangle. 
\end{align*}
Moreover, we have 
\[\langle Q_p(v_2),\tilde{v}_2 \rangle = \langle Q_p(v_1),\tilde{v}_1 \rangle = \varphi(p,f(p)) \geq \frac{\theta^2}{L}\] 
by (\ref{obliqueness}). Putting these facts together, we conclude that $|b| \leq C$ for some uniform constant $C$.

\textit{Step 5:} In the last step, we show that $|Df_p| \leq C$ for all points $p \in \Omega$. As above, we define a function $\beta: \Sigma \to \mathbb{R}$ by 
\[\beta(p,f(p)) = \frac{2}{\sqrt{\det(I + Df_p^* \, Df_p)}}.\] 
It follows from Lemma \ref{pde.for.beta} that the function $\beta$ satisfies the inequality 
\[\Delta_\Sigma \beta \leq -\kappa \, \beta \, (1 - \beta^2).\] 
This gives 
\begin{equation} 
\label{pde.for.log.beta}
\Delta_\Sigma (\log \beta) \leq -\kappa \, (1 - \beta^2). 
\end{equation}
Moreover, the restriction $\beta|_{\partial \Sigma}$ is uniformly bounded from below. Using (\ref{pde.for.log.beta}) and the maximum principle, one obtains a uniform lower bound for $\inf_{p \in \Omega} \beta(p,f(p))$. This completes the proof of Proposition \ref{estimate.for.Df}. \\

After these preparations, we now sketch the proof of Theorem \ref{existence.uniqueness}. The proof uses the continuity method. We first construct domains $\Omega_t,\tilde{\Omega}_t \subset N$ with the following properties: 
\begin{itemize}
\item For each $t \in (0,1]$, the domains $\Omega_t$ and $\tilde{\Omega}_t$ are uniformly convex, and $\text{\rm area}(\Omega_t) = \text{\rm area}(\tilde{\Omega}_t)$.
\item $\Omega_1 = \Omega$ and $\tilde{\Omega}_1 = \tilde{\Omega}$.
\item If $t \in (0,1]$ is sufficiently small, then $\Omega_t$ and $\tilde{\Omega}_t$ are geodesic disks in $N$. Moreover, the radius converges to $0$ as $t \to 0$.
\end{itemize}
In order to construct domains $\Omega_t,\tilde{\Omega}_t \subset N$ with these properties, we consider the sub-level sets of suitable boundary defining functions (see \cite{Brendle} for details). We then consider the following problem: \\

$(\star_t)$ \textit{Find all area-preserving minimal maps $f: \Omega_t \to \tilde{\Omega}_t$ that map a given point on the boundary of $\Omega_t$ to a given point on the boundary of $\tilde{\Omega}_t$.} \\

As $t \to 0$, the domains $\Omega_t$ and $\tilde{\Omega}_t$ converge to the unit disk $\mathbb{B}^2 \subset \mathbb{R}^2$ after rescaling. Hence, for $t \to 0$, the problem $(\star_t)$ reduces to the problem of finding all area-preserving minimal maps from the flat unit disk $\mathbb{B}^2$ to itself. This problem is well understood: in fact, an area-preserving map from $\mathbb{B}^2$ to itself is minimal if and only if it is a rotation. 

Using Proposition \ref{estimate.for.Df}, we obtain uniform a priori estimates for solutions of $(\star_t)$. Moreover, it turns out that each solution of $(\star_t)$ is non-degenerate in the sense that the linearized operator is invertible. Hence, it follows from standard continuity arguments that $(\star_t)$ has a unique solution for each $t \in (0,1]$.

\section{The Lagrangian mean curvature flow}

In this final section, we briefly discuss the flow approach to special Lagrangian geometry. To that end, we consider a Lagrangian submanifold of a K\"ahler manifold $(M,g)$, and evolve it by the mean curvature flow. It was shown by Smoczyk that a Lagrangian submanifold of a K\"ahler-Einstein manifold remains Lagrangian when evolved by the mean curvature flow:

\begin{theorem}[K.~Smoczyk \cite{Smoczyk1},\cite{Smoczyk2}]
Let $(M,g)$ be a K\"ahler-Einstein manifold, and let $\{\Sigma_t: t \in [0,T)\}$ be a family of closed submanifolds of $(M,g)$ which evolve by the mean curvature flow. If $\Sigma_0$ is Lagrangian, then $\Sigma_t$ is Lagrangian for all $t \in [0,T)$.
\end{theorem}

It is a very interesting question to study the longtime behavior of the Lagrangian mean curvature flow. Thomas and Yau \cite{Thomas-Yau} conjectured that the flow exists for all time provided that the initial surface $\Sigma_0$ satisfies a certain stability condition. Examples of finite-time singularities were recently constructed by Neves \cite{Neves}.

In the following, we discuss some results about Lagrangian graphs evolving by mean curvature flow. The case of graphs is much better understood than the general case, and some strong results are known in this setting. Let us first consider the torus $\mathbb{T}^{2n} = \mathbb{R}^{2n} / \mathbb{Z}^{2n}$. We assume that $\mathbb{R}^{2n}$ is equipped with its standard metric and complex structure, so that $J \frac{\partial}{\partial x_k} = \frac{\partial}{\partial y_k}$ and $J \frac{\partial}{\partial y_k} = -\frac{\partial}{\partial x_k}$. The torus $\mathbb{T}^{2n}$ inherits a metric and complex structure in the standard way. We then consider submanifolds of the form 
\[\Sigma = \{(p,f(p)): p \in \mathbb{T}^n\},\] 
where $f$ is a smooth map from $\mathbb{T}^n$ to itself. The submanifold $\Sigma$ is Lagrangian if and only if the map $f$ can locally be written in the form $f = \nabla u$ for some potential function $u$. Smoczyk and Wang were able to analzye the longtime behavior of the mean curvature flow in the special case when the potential function $u$ is convex.

\begin{theorem}[K.~Smoczyk, M.T.~Wang \cite{Smoczyk-Wang}]
\label{flat.case}
Let $\Sigma_0$ be a Lagrangian submanifold of $\mathbb{T}^{2n}$ which can be written as the graph of a map $f_0: \mathbb{T}^n \to \mathbb{T}^n$. Moreover, suppose that the eigenvalues of $(Df_0)_p$ are strictly positive for each point $p \in \mathbb{T}^n$. Finally, let $\{\Sigma_t: t \in [0,T)\}$ denote the unique maximal solution of the mean curvature flow with initial surface $\Sigma_0$. Then $T = \infty$, and the surfaces $\Sigma_t$ converge to a totally geodesic Lagrangian submanifold as $t \to \infty$.
\end{theorem}

We next consider the Lagrangian mean curvature flow in a product manifold.

\begin{theorem}[M.T.~Wang \cite{Wang2}]
\label{riemann.surface}
Let $N$ and $\tilde{N}$ be compact Riemann surfaces with the same constant curvature $c$. Moreover, suppose that $f_0: N \to \tilde{N}$ is an area-preserving diffeomorphism, and let 
\[\Sigma_0 = \{(p,f_0(p)): p \in N\} \subset N \times \tilde{N}\] 
denote the graph of $f_0$. Finally, let $\{\Sigma_t: t \in [0,T)\}$ be the unique maximal solution of the mean curvature flow with initial surface $\Sigma_0$. Then $T=\infty$, and each surface $\Sigma_t$ is the graph of an area-preserving diffeomorphism $f_t: N \to \tilde{N}$. Finally, the maps $f_t$ converge smoothly to an area-preserving minimal map as $t \to \infty$.
\end{theorem}

The same result was proved independently by Smoczyk \cite{Smoczyk2} under an extra condition on the Lagrangian angle.

Theorem \ref{riemann.surface} gives a new proof of the existence of minimal maps between Riemann surfaces; the existence of such maps was established earlier by Schoen \cite{Schoen} using harmonic map techniques. A stronger result holds when $N = \tilde{N} = S^2$: 

\begin{theorem}[M.T.~Wang \cite{Wang2}]
\label{sphere}
Let $f_0$ be an area-preserving diffeomorphism from $S^2$ to itself, and let 
\[\Sigma_0 = \{(p,f_0(p)): p \in S^2\} \subset S^2 \times S^2\] 
denote the graph of $f_0$. Moreover, let $\{\Sigma_t: t \in [0,T)\}$ be the unique maximal solution of the mean curvature flow with initial surface $\Sigma_0$. Then $T=\infty$, and each surface $\Sigma_t$ is the graph of an area-preserving diffeomorphism $f_t: S^2 \to S^2$. Finally, the maps $f_t$ converge to an isometry of $S^2$ as $t \to \infty$.
\end{theorem}

The proofs of Theorems \ref{flat.case} -- \ref{sphere} rely on maximum principle arguments. These techniques also have important applications to the study of area-decreasing maps between spheres (cf. \cite{Tsui-Wang2}, \cite{Wang1}). A detailed discussion of the Lagrangian mean curvature flow can be found in \cite{Wang-survey}.

In a remarkable paper, Medo\v s and Wang \cite{Medos-Wang} generalized this result to higher dimensions. In higher dimensions, it is necessary to impose a pinching condition on the initial map $f_0$:

\begin{theorem}[I.~Medo\v s, M.T.~Wang \cite{Medos-Wang}]
\label{complex.projective.space}
Given any positive integer $n$, there exists a real number $\Lambda(n) > 1$ such that the following holds: Let $f_0: \mathbb{CP}^n \to \mathbb{CP}^n$ be a symplectomorphism satisfying 
\[\frac{1}{\Lambda(n)} \, |v| \leq |Df_p(v)| \leq \Lambda(n) \, |v|\] 
for all vectors $v \in T_p \mathbb{CP}^n$. Moreover, let 
\[\Sigma_0 = \{(p,f(p)): p \in \mathbb{CP}^n\} \subset \mathbb{CP}^n \times \mathbb{CP}^n\] 
denote the graph of $f_0$, and let $\{\Sigma_t: t \in [0,T)\}$ be the unique maximal solution of the mean curvature flow with initial surface $\Sigma_0$. Then $T = \infty$, and each surface $\Sigma_t$ is the graph of a symplectomorphism $f_t: \mathbb{CP}^n \to \mathbb{CP}^n$. Moreover, the maps $f_t$ converge smoothly to a biholomorphic isometry of $\mathbb{CP}^n$ as $t \to \infty$.
\end{theorem}

In the remainder of this section, we sketch the main ingredients involved in the proof of Theorem \ref{complex.projective.space} (see \cite{Medos-Wang} for details). For each $t \geq 0$, one defines a function $\beta_t: \Sigma_t \to \mathbb{R}$ by 
\[\beta_t = \prod_{k=1}^{2n} \frac{1}{\sqrt{1+\lambda_k^2}},\] 
where $\lambda_1,\hdots,\lambda_n$ denote the singular values of $Df_t$. Since $f_t$ is a symplectomorphism, the singular values of $Df_t$ occur in pairs of reciprocal numbers. We may therefore assume that $\lambda_i \lambda_{\tilde{i}} = 1$, where $\tilde{i} = i + (-1)^{i-1}$. Consequently, $\beta_t \leq 2^{-n}$, and equality holds if and only if $\lambda_1 = \hdots = \lambda_n = 1$.

The function $\beta_t$ satisfies an evolution equation of the form 
\begin{align*} 
\frac{\partial}{\partial t} \beta_t 
&= \Delta_{\Sigma_t} \beta_t + \frac{\beta_t}{2} \sum_{k=1}^{2n} \Big ( \frac{1-\lambda_k^2}{1+\lambda_k^2} \Big )^2 \\ 
&+ \beta_t \sum_{i,j,k=1}^{2n} h_{ijk}^2 - 2\beta_t \sum_{k=1}^{2n} \sum_{i<j} (-1)^{i+j} \, \lambda_i \, \lambda_j \, (h_{i\tilde{i}k} \, h_{j\tilde{j}k} - h_{i\tilde{j}k} \, h_{j\tilde{i}k}) 
\end{align*}
where $h_{ijk} = \langle I\!I(e_i,e_j),Je_k \rangle$ denote the components of the second fundamental form of $\Sigma_t$ (cf. \cite{Medos-Wang}, Proposition 2). In order to apply the maximum principle to the function $\beta_t$, one needs to verify that the terms on the right hand side of the evolution equation are nonnegative. In fact, it is shown in \cite{Medos-Wang} that 
\begin{equation} 
\label{2nd.fundamental.form}
\sum_{i,j,k=1}^{2n} h_{ijk}^2 - 2 \sum_{k=1}^{2n} \sum_{i<j} (-1)^{i+j} \, \lambda_i \, \lambda_j \, (h_{i\tilde{i}k} \, h_{j\tilde{j}k} - h_{i\tilde{j}k} \, h_{j\tilde{i}k}) \geq \delta \sum_{i,j,k=1}^{2n} h_{ijk}^2, 
\end{equation}
provided that the singular values $\lambda_1,\hdots,\lambda_n$ are sufficiently close to $1$. In order to verify this, Medo\v s and Wang consider the quadratic form 
\[\mathcal{Q}(h) = \sum_{i,j,k=1}^{2n} h_{ijk}^2 - 2 \sum_{k=1}^{2n} \sum_{i<j} (-1)^{i+j} \, (h_{i\tilde{i}k} \, h_{j\tilde{j}k} - h_{i\tilde{j}k} \, h_{j\tilde{i}k}).\] 
The estimate (\ref{2nd.fundamental.form}) is then a consequence of the following result (cf. \cite{Medos-Wang}, Lemma 4): 

\begin{proposition}
The quadratic form $\mathcal{Q}(h)$ satisfies 
\begin{equation} 
\label{estimate.for.Q}
\mathcal{Q}(h) \geq \frac{2}{9} \sum_{i,j,k=1}^{2n} h_{ijk}^2. 
\end{equation}
\end{proposition}

In order to prove the inequality (\ref{estimate.for.Q}), we observe that $\sum_{i=1}^{2n} (-1)^i \, h_{i\tilde{i}k} = 0$ for each $k$. From this, we deduce that $\sum_{i,j=1}^{2n} (-1)^{i+j} \, h_{i\tilde{i}k} \, h_{j\tilde{j}k} = 0$ for each $k$. Consequently, the quadratic form $\mathcal{Q}(h)$ can be rewritten as 
\begin{align*} 
\mathcal{Q}(h) 
&= \sum_{i,j,k=1}^{2n} h_{ijk}^2 - \sum_{i,j,k=1}^{2n} (-1)^{i+j} \, (h_{i\tilde{i}k} \, h_{j\tilde{j}k} - h_{i\tilde{j}k} \, h_{j\tilde{i}k}) \\ 
&= \sum_{i,j,k=1}^{2n} h_{ijk}^2 + \sum_{i,j,k=1}^{2n} (-1)^{i+j} \, h_{i\tilde{j}k} \, h_{j\tilde{i}k} \\ 
&= \frac{1}{2} \sum_{i,j,k=1}^{2n} \big ( (-1)^i \, h_{i\tilde{j}k} + (-1)^j h_{\tilde{i}jk} \big )^2. 
\end{align*}
On the other hand, the identity 
\begin{align*} 
2 \, h_{ijk} 
&= (-1)^i \, \big ( (-1)^i \, h_{ijk} + (-1)^{\tilde{j}} \, h_{\tilde{i}\tilde{j}k} \big ) \\ 
&+ (-1)^i \, \big ( (-1)^i \, h_{ijk} + (-1)^{\tilde{k}} \, h_{\tilde{i}j\tilde{k}} \big ) \\ 
&+ (-1)^{i+j+k} \, \big ( (-1)^k \, h_{\tilde{i}\tilde{j}k} + (-1)^j \, h_{\tilde{i}j\tilde{k}} \big ) 
\end{align*}
implies 
\begin{align*} 
4 \, h_{ijk}^2 
&\leq 3 \, \big ( (-1)^i \, h_{ijk} + (-1)^{\tilde{j}} \, h_{\tilde{i}\tilde{j}k} \big )^2 \\ 
&+ 3 \, \big ( (-1)^i \, h_{ijk} + (-1)^{\tilde{k}} \, h_{\tilde{i}j\tilde{k}} \big )^2 \\ 
&+ 3 \, \big ( (-1)^k \, h_{\tilde{i}\tilde{j}k} + (-1)^j \, h_{\tilde{i}j\tilde{k}} \big )^2. 
\end{align*} 
Summation over $i,j,k$ yields 
\[4 \sum_{i,j,k=1}^{2n} h_{ijk}^2 \leq 18 \, \mathcal{Q}(h),\] 
as claimed.


\begin{thebibliography}{99}
\bibitem{Brendle}
S.~Brendle, \textit{Minimal Lagrangian diffeomorphisms between domains in the hyperbolic plane,} J. Diff. Geom. 80, 1--22 (2008)

\bibitem{Brendle-Warren}
S.~Brendle and M.~Warren, \textit{A boundary value problem for minimal Lagrangian graphs,} J. Diff. Geom. 84, 267--287 (2010)

\bibitem{Caffarelli}
L.~Caffarelli, \textit{Boundary regularity of maps with convex potentials, II,} Ann. of Math. 144, 453--496 (1996)

\bibitem{Caffarelli-Nirenberg-Spruck}
L.~Caffarelli, L.~Nirenberg, and J.~Spruck, \textit{The Dirichlet problem for nonlinear second order elliptic equations, III: functions of the eigenvalues of the Hessian,} Acta Math. 155, 261--301 (1985)

\bibitem{Delanoe}
P.~Delano\"e, \textit{Classical solvability in dimension two of the second boundary-value problem associated with the Monge-Amp\`ere operator,} Ann. Inst. H. Poincar\'e 8, 443--457 (1991)

\bibitem{Harvey-Lawson}
R.~Harvey and H.B.~Lawson, Jr., \textit{Calibrated geometries,} Acta Math. 148, 47--157 (1982)

\bibitem{Medos-Wang}
I.~Medo\v s and M.T.~Wang, \textit{Deforming symplectomorphisms of complex projective spaces by the mean curvature flow,} J. Diff. Geom. 87, 309--342 (2011)

\bibitem{Neves}
A.~Neves, \textit{Finite time singularities for Lagrangian mean curvature flow,} arxiv:1009.1083

\bibitem{ONeill}
B.~O'Neill, \textit{Semi-Riemannian geometry,} Academic Press, New York (1983)

\bibitem{Schoen}
R.~Schoen, \textit{The role of harmonic mappings in rigidity and deformation problems,} Complex geometry, Proc. Osaka International Conference, Marcel Dekker, New York, 1993

\bibitem{Smoczyk1}
K.~Smoczyk, \textit{Der Lagrangesche mittlere Kr\"ummungsflu\ss{},} Habilitationsschrift, Leipzig University (1999)

\bibitem{Smoczyk2}
K.~Smoczyk, \textit{Angle theorems for the Lagrangian mean curvature flow,} Math. Z. 240, 849--883 (2002)

\bibitem{Smoczyk-Wang}
K.~Smoczyk and M.T.~Wang, \textit{Mean curvature flows of Lagrangian submanifolds with convex potentials,} J. Diff. Geom. 62, 243--257 (2002)

\bibitem{Thomas-Yau}
R.P.~Thomas and S.T.~Yau, \textit{Special Lagrangians, stable bundles, and mean curvature flow,} Comm. Anal. Geom. 10, 1075--1113 (2002)

\bibitem{Tsui-Wang1}
M.P.~Tsui and M.T.~Wang, \textit{A Bernstein type result for special Lagrangian submanifolds,} Math. Res. Lett. 9 529--535 (2002)

\bibitem{Tsui-Wang2}
M.P.~Tsui and M.T.~Wang, \textit{Mean curvature flows and isotopy of maps between spheres,} Comm. Pure Appl. Math. 57, 1110--1126 (2004)

\bibitem{Urbas1}
J.~Urbas, \textit{On the second boundary value problem for equations of Monge-Amp\`ere type,} J. Reine Angew. Math. 487, 115--124 (1997)

\bibitem{Urbas2}
J.~Urbas, \textit{The second boundary value problem for a class of Hessian equations,} Comm. PDE 26, 859--882 (2001)

\bibitem{Urbas3}
J.~Urbas, \textit{A remark on minimal Lagrangian diffeomorphisms and the Monge-Amp\`ere equation,} Bull. Austral. Math. Soc. 76, 215--218 (2007)

\bibitem{Wang1}
M.T.~Wang, \textit{Mean curvature flow of surfaces in Einstein four-manifolds,} J. Diff. Geom. 57, 301--338 (2001)

\bibitem{Wang2}
M.T.~Wang, \textit{Deforming area-preserving diffeomorphisms of surfaces by mean curvature flow,} Math. Res. Letters 8, 651--662 (2001)

\bibitem{Wang-survey}
M.T.~Wang, \textit{Some recent developments in Lagrangian mean curvature flows,} Surveys in Differential Geometry vol. XII, pp.~333--347, International Press, Somerville MA (2008)

\bibitem{Yuan}
Y.~Yuan, \textit{A Bernstein problem for special Lagrangian equations,} Invent. Math. 150, 117--125 (2002)
\end{thebibliography}
\end{document}